# Complete invariant graphs of alternating knots

Christian Soulié

First submission: April 2004 (revision 1)

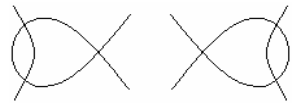

**Abstract :**    Chord diagrams and related enlacement graphs of alternating knots are enhanced to obtain  complete invariant graphs including chirality detection. Moreover, the equivalence by common enlacement graph is specified and the neighborhood graph is defined for general purpose and for special application to the knots.

### I - Introduction :

Chord diagrams are enhanced to integrate the state sum of all flype moves and then produce an invariant graph for alternating knots. By adding local writhe attribute to these graphs, chiral types of knots are distinguished. The resulting chord-weighted graph is a complete invariant of alternating knots. Condensed chord diagrams and condensed enlacement graphs are introduced and a new type of graph of general purpose is defined : the neighborhood graph. The enlacement graph is enriched by local writhe and chord orientation. Hence this enhanced graph distinguishes mutant alternating knots. As invariant by flype it is also invariant for all alternating knots. The equivalence class of knots with the same enlacement graph is fully specified and extended mutation with flype of tangles is defined. On this way, two enhanced graphs are proposed as complete invariants of alternating knots.

## II - Definitions and condensed graphs :

### II-1 Knots

A *lacet* is a closed curve embedded on a 2-manifold. A *real knot* is a closed curve embedded on R^3 or S^3 without intersecting itself, up to ambient isotopy. The embedding type of a set of several closed curves into a 3-manifold without any intersection between them and themselves is a *link*.

The projection of a knot or a link onto a 2-manifold is considered with all curves in general position and all multiple points are transversal double points which will be called *crossing points*. Such a projection is called the *shadow* by the knots theorists [Ad],[So], or the *lacet* by the graph theorists [CrRo] . By adding a 3D-attribute to all crossing points of a lacet to distinguish the over strand from the under strand, we obtain the *knot diagram*. The lacet is a 4-regular plane graph whose vertices are the crossing points and whose edges are the strands between 2 crossing points. For a knot K we will note this graph drawing G(K). For more details see [Ad],[So],[Rob]. In the sequel we'll consider knot diagrams on S² surface only.

When a topological sphere is cut by the knot in 2 points, it separates it into 2 different parts called *1-tangles* or components of the knot. When the sphere cuts in n points, it separates the knot into two *n-tangles*. When a knot cannot be split into any 1-tangles, it is *prime*. Horst Schubert demonstrated in 1949 that there exists for any knot a unique composition into prime knots.

The path on G(K) which is transversal to each vertex is generating a word W(K). When the letters of this word are the labels of the crossing points with (resp. without) attached 3D-attributes (knot or resp. lacet), then it is called the *Gauss code* [Ga]. A graph is said *Gaussian* if there is an embedding which contains a transversal Eulerian circuit [KeSe]. Two different directions can be chosen for this path. An oriented knot is a knot where one direction of the path is fixed. The change of this orientation is called *reversion*. When the sequence of over/under 3D-attributes along the path is alternating the knot is *alternating*. It is well known that an alternating knot is reduced when removing all the loops (successive occurrence of the same letter in Gauss code). Our study is mainly limited to reduced alternating knots in the sequel, without to lose on generality for complete invariants of isotopic knots.

The knot diagram obtained by switching all 3D-attributes is the diagram of the mirror image of the knot (a view from behind the mirror). There are exactly 2 different knots K and its mirror image K* corresponding to one lacet : G(K)=G(K*). The theorem of Gauss [Ga], reported later here as a condition of realizability of W(K) says that every proper subcircuit of the Gaussian (Eulerian) circuit is of odd length. A consequence of this theorem is that alternating attributes are always possible to apply on vertices [KeSe].

When the knot is invariant by mirroring it is *achiral* : K=K*. When mirror image is different (cannot be obtained by ambient isotopy) from its preimage the knot is chiral and we must distinguish K and K*:   K≠K*. In the sequel we will not distinguish achiral knots where chirality depends on orientation and oriented non-invertible knots. Orientation is necessary to identify global properties where local geometric figures are not enough but orientation is not part of the nature of a knot unless by adding it to the definition.

The number of crossing points is c and a knot is written $C_n$ where n is the number given by Rolfsen in his classification table or, for higher c values, by Thistlesthwaite table. See figure 1 for an example of the chiral knot $7_7$.

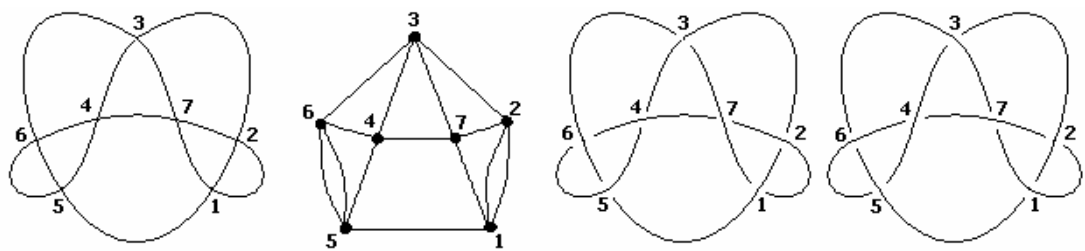

Figure 1 . G($7_7$) lacet and graph, knots diagrams of $7_7$ and $7_7$*.

**II-2 Sign of crossing points**

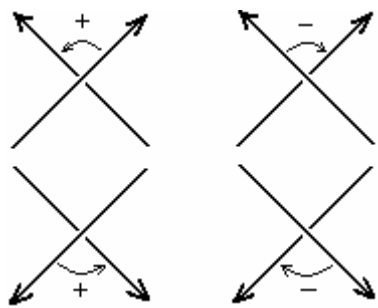

Figure 2 . Sign of crossing

We say that a crossing is positive (resp. negative) on an oriented knot diagram as shown on figure 2. There is a value of this definition for non oriented knots because the sign of the crossing doesn't change when reversing the orientation. We can then define a 3D attribute ε and ε =+1 (resp. ε =-1) if the crossing is positive (resp. negative). It is well known that the sum of this attribute over all internal crossing is the *writhe* invariant and the sum of ε attribute over all crossings between 2 components of a link, divided by 2, is also an invariant: the *linking number*. That is the reason why ε is also called the *local writhe*.

Nevertheless the author notices that not only the sum but also every local writhe is an invariant for any move of the knot diagram which corresponds to an ambient isotopy, like for example: Reidemeister of type III (R3) and Flype (F). This can be done by an adequate identification (labels) of crossings before and after the move and it is well defined by $C^n$ ($n \geq 2$) continuities of the curve and of the move.

The invariant nature of the sign of the crossing on oriented 2-manifolds is due to the fact that it gives the coorientation between the crossing and the vector N topologically normal to the surface and pointing to the point of view. If O (resp. U) is the vector tangent to the over (resp. under) strand getting for direction the orientation of the knot, then (O,U,N) is a direct (resp. inverse) trihedron when the sign of crossing is positive (resp. negative).

**II-3 Chord diagrams**

After pointing the diagram of the lacet by an origin and choosing an orientation, we walk on a path along the curve until returning back to the origin and then generate a word W(K) which is the sequence of the crossings in the order we meet them on the path. W(K) is a double occurrence word. This word is built with the labels of the c crossings as letters of the alphabet. The lacet is not changed by circular permutation of the letters of its word: it's like changing the choice of the interval between two crossings where the root is placed. The lacet is unchanged when we maintain the coorientation of the lacet with the surface supporting it but we obtain the mirror image of it if we reverse this coorientation. W(K) is the *Gauss code* of the lacet and of the knot diagram with the same shadow.

When adding 3D attributes over (O) and under (U) to each crossing letter depending of how it occurs in the sequence, we get a word W'(K) of the knot K. Like its shadow (the lacet), the knot is not changed by circular permutation nor by reversing coorientation. W(K) = W(K*) and W'(K) = W'(K*).

On example $7_7$ of figure 1 and 3, we get :
W($7_7$) = 1 2 3 4 5 6 4 7 2 1 7 3 6 5
W'($7_7$) = O1 U2 O3 U4 O5 U6 O4 U7 O2 U1 O7 U3 O6 U5
W'($7_7$*) = U1 O2 U3 O4 U5 O6 U4 O7 U2 O1 U7 O3 U6 O5

W(K) is *Gauss code* of the lacet and of the knot diagram with the same shadow.
W'(K) is *Gauss code* of the knot diagram and then one of Gauss codes of a knot.

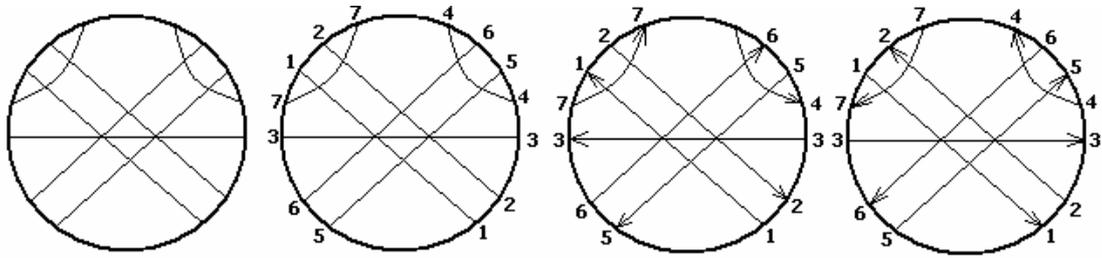

Figure 3 . Chord diagrams of $7_7$ lacet without and with labeling, of $7_7$ knot with its mirror image.

If we put the labels of the crossings on a circle in the order of the Gauss code and if we join by a chord all pairs of identical labels then we obtain a labeled *chord diagram* as on the second diagram of figure 3. The first diagram is unlabeled and then is a better representation of the lacet (class under equivalence of letters permutations). In short, we will write CD(K) for it. Chord diagrams are a good representation of the relations of intersections between geometric or topologic objects. It can be used for braids, for Jordan curve segments, for connectors... In the particular case of lacets or knots we call them *Gauss diagrams,* whenever confusion is possible. A CD is a 3-regular simple connected Hamiltonian graph.

3D attributes are brought on the chord diagram by drawing arrows instead of chords. By convention the tail of an arrow is attached to the over part of the strand and the head to the under part. On figure 3 we see that $7_7$ lacet is achiral where we know that $7_7$ is chiral by the effect of the 3D attributes. Chirality is then related to both lacet geometry and combinatorial arrangement of attributes. The reader will notice that CDs with arrow-chords doesn't distinguish $7_7$ from $7_7$*. Effectively, one CD is the same of the other when seen from behind. CD is a 3-regular abstract graph: it is not depending of its embedding on the plane of drawing. When rooted and oriented, CD is bijective with Gauss code which defines only adjacencies of CD graph. The place of the root corresponds to the circular permutation of the word, and orientation reversion to word inversion.

If we use the signs of crossings, we can assign them to the chords or duplicate to both ends of each chord. We can also put sign attributes in the Gauss code:
W"($7_7$) = +1+2-3-4+5+6-4-7+2+1-7-3+6+5
W"($7_7$*) = -1-2+3+4-5-6+4+7-2-1+7+3-6-5

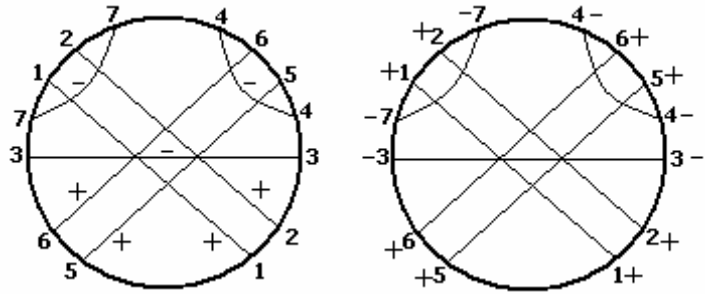

Figure 4 . Two notations of sign attributes on CD($7_7$)

When projecting a knot onto a sphere instead of the plane we get far less diagrams. Every knot diagram on the plane is got from the one of the sphere by selecting a face (a disc on $S^2 \setminus K$) of G(K) and use it for the outer face of the drawing on the plane. It is equivalent to select a puncture point on the sphere. Some people, like Goussarov, are selecting the puncture point on a segment of the knot [GPV]. A *long knot* is obtained by this way with one point to infinity. The circle of the CD is then open and we can represent the long knot with a linear chord diagram (LCD), image of an open knot or a self-intersecting Jordan curve in 3-space.

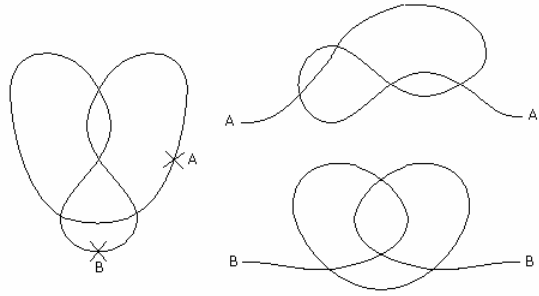

Figure 5 . Open knots or long knots of the figure eight knot

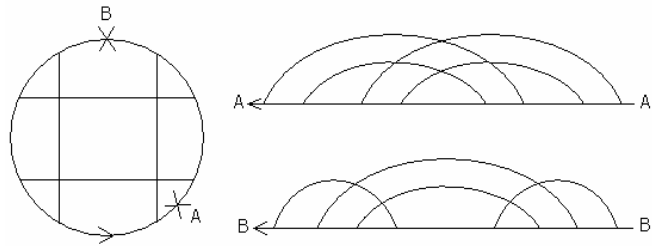

Figure 6 . Linear Chord Diagrams (LCD) of figure8-knot for cuts A and B
For the long knots : A=∞ and B=∞

If we call *S-isotopy* the spherical isotopy (ambient isotopy) between 2 plane knot diagrams which leaves spherical knot diagram the same, we see that a CD doesn't distinguish the outer face and then it is remarkable to notice that a chord diagram is invariant by S-isotopy.
Proof :   The cyclic ordering at any vertex of G(K) is not modified when changing the outer face of the plane knot diagram. The successor of an edge in a Gauss path is the adjacent one obtained by traversal. The successor by  traversal is unchanged if the cyclic ordering of adjacent vertices of the crossing is unchanged and that is the case.   ÿ
This property is the same for the CD of a link as the property of path invariance (cyclic ordering) is the same for a crossing between components (*inter*) and for a self crossing (*intra*). Once the choice of left or right turn is chosen for each inter-crossing, this path is invariant by S-isotopy for each vertex. By this construction the path is a Gaussian circuit and defines a "single" CD. Choice of the position of the root of each component and relative orientation of each component is giving a different CD. That is the reason why we limit our study to knots (links with one component). For definition of CDs of a link, see [Ka], [ReRo].

**Theorem 1-1** : *The chord diagram of a lacet, a knot or a link is invariant by S-isotopy.*

**II-4 Enlacement graphs**

Let us say that two letters, x and y, are *enlaced* in a 2-occurrence word like the Gauss code if they occur in the following order :   ...x...y...x...y.... There is one occurrence of a letter between two occurrences of the other, when these letters are enlaced. The *enlacement graph* has one vertex for each letter of a 2-occurrence word and an edge between 2 vertices when the corresponding letters are enlaced. We'll use in the sequel LG(K) for the enlacement graph of a knot K. For LG(K), the letters are the labels of the crossings in the Gauss code W(K) and its chord diagram CD(K). From each crossing on the knot diagram G(K) are issued 2 loops, each with a sequence of crossings. The crossings which appear one time in a sequence of a loop are enlaced with the crossing generating these loops. This relation of enlacement is the same for the 2 loops and then is commutative. The initial name given by the author was "Loop Graph" (LG). We must notice here that two chords are intersecting on a CD when the related crossings are enlaced. Then an equivalent definition of LG(K) is that it is the chord-intersection graph of CD(K). For the same reason LG(K) is also the chord-intersection graph of LCD(K).

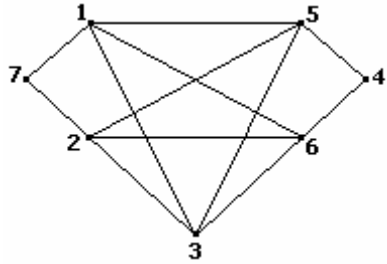

Figure 7 . LG ($7_7$)

We give here the definition of the chord-intersection graph as written by Blake Mellor [BMe]. It is generalizing the previous definition because it is defined as a map from the set of all CDs to the set of simple graphs. Not all the CDs can be realized as a lacet on the sphere but we may generalize as Louis Kauffman did for knot diagrams [Ka] with virtual lacets, covering by this way all types of CDs. Of course the set of CDs exist by itself (i.e. out of knot theory and its representation by a knot, a link, or a weight system in Vassiliev-Goussarov-Kontsevich theory) and then the related set of LGs as well.

**Definition :** *Given a chord diagram CD, we define its intersection graph LG(CD) as the graph such that :*
- *LG(CD) has a vertex for each chord of CD*
- *Two vertices of LG(CD) are connected by an edge if and only if the corresponding chords in CD intersect, i.e. their endpoints on the bounding circle alternate.*

The chord-intersection graph is an abstract graph.
This intersection graph can be enhanced by chord attributes like local writhe on two different ways : simply by attaching these attributes to the related vertices, or by assigning some product of these attributes to the edges.
When the chords of a CD are oriented, we can generate from this a directed LG (counterclockwise direction of rotation from one arrow head to the other in the CD, when arrow-chords are intersecting).

**II-5 Condensed graphs**

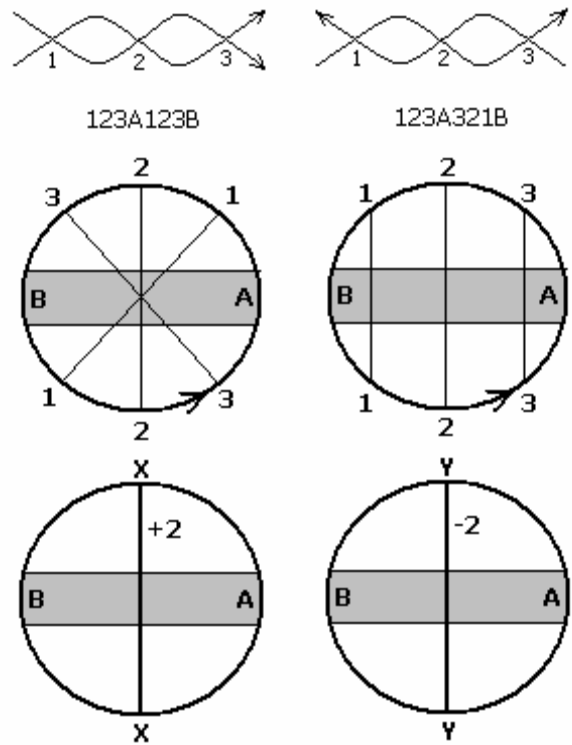

Figure 7 . Torsads with related secant and parallel sets of chords.
Corresponding condensed chord diagrams with weighted edges.

We call a *torsad* a sequence of half twists. A *n-torsad* has n half twists. A torsad n-$t^k$ is a regular braid (n twists on a trivial braid) with k strings. We call $t^2$ a *digonal torsad* (two strings) and $t^2_2$ a *digon* (digonal torsad with two crossings). So is $t^3$ a trigonal torsad and 1-$t^3$ or $t^3_3$ is a *trigon* which surround a triangle. A torsad n-$t^k$ is a particular k-tangle with k.n crossings. Then we have equality of the notations n-$t^k = t^k_{n.k}$ . We may generalize the definition by using fractional twists instead of half twists. Then a (p/k)-$t^k$ torsad is a sequence of p/k half twists which is endowed of p crossings. Then, (p/k)-$t^k = t^k_p$ . $t^2$ digonal torsads are the elements of Conway construction of rational tangles and links.

We define a torsad $t^2$ to be *secant* when the strings are oriented in the same direction. This property is well defined as it doesn't depend on the orientation of the knot. That is not the case for a link whose component orientations are independant and of arbitrary choice. The chords representing the crossings internal to a secant torsad are intersecting on CD. See figure 7. On an analog way, a torsad $t^2$ is *parallel* when the strings are oriented in opposite directions and the chords of its crossings are parallel on CD. We prefer "secant" and "parallel" to "positive" and "negative" as commonly used in the literature to avoid any confusion with the sign (local writhe) of a crossing. Indeed a 1-$t^2$ or $t^2_1$ torsad is a simple crossing which can be positive (local writhe) or negative. We can extend local writhe to the whole torsad when it is alternating because in that case, the crossing sign is constant. Note that it is the same for a non-alternating reduced link. Otherwise the link could be reduced by a Reidemeister move of type II. There are no non-trivial moves inside alternating torsads on a link diagram. So we can define positive secant, negative secant, positive parallel and negative parallel $t^2$ torsads. An attribute like positive/secant is mainly/purely global : it depends partly/only (local writhe is a combination of local over/under attribute and global : codirection of the set of strings at the crossing) on the complementary 2-tangle connector class (see chapter IV-5 the classification by connectors). If we don't mention it in the sequel we'll consider each torsad to be digonal.

Each bundle of chords in a CD can be replaced or packed ("condensed") into a weighted chord. The weight of a n-torsad is +(n-1) when secant, -(n-1) when parallel, 0 when it is a single crossing (n=1). A single crossing is both secant and parallel as 0 is both positive and negative. The same condensation can be applied to the enlacement graph. See on figure 8 below. You may compare with Conway notation where the condensed LG adds a graph structure to Conway coefficients and is a more homogeneous notation with non rational and non algebraic knots. For example, compare condensed LG($7_7$) with 21112 Conway notation.

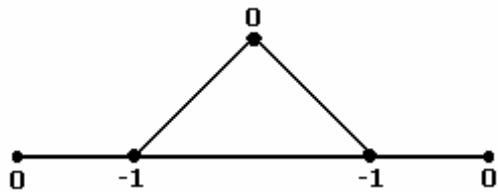

Figure 8 . Condensed LG($7_7$)

A condensed LG is a vertex-weighted simple graph (no loops and no multiple edges). A positive vertex of the condensed LG is a clique (complete subgraph $K_n$) in the original LG, a negative vertex is a disjoint set of vertices (an "*anticlique*"), and an edge is a bipartite subgraph joining cliques and anticliques.

Here we used a new type of graph which could be of general purpose in the graph theory: it is the graph of common neighborhood. We'll call it *Neighborhood Graph N(G)* and it is obtained from any graph G by identifying vertices with common neighborhood and removing resulting loops and duplicate edges. We assign an integer weight to the vertices and give a zero weight for each vertex of G. Then, each time we identify 2 adjacent vertices, the weight of the contracted vertex is the sum of the respective (positive) weight plus one. If we contract 2 disjoint vertices, the weight of the contracted vertex is the sum of the (negative) weights minus one. There is a unique end to this process where the graph cannot be contracted anymore. The irreducible contracted graph is N(G). N is a projective map from the set of all graphs to some subset of integer-vertex-weighted graphs. Chromatic number, clique number, clique cover number, and stability number of G can be extracted easily from the neighborhood graph N (G).

Instead of contracting edges without any condition like for minors, we contract edges and non-edges under the condition of common neighborhood between vertices. This gives a lot more symmetry between edges and non-edges in the structure of N(G), also between cliques and anticliques, and the complement of the neighborhood graph is the neighborhood graph of the complement. Common neighborhood of subsets is detected in a graph when a subgraph is connected to the remaining part of the graph only by one bipartite complete subgraph. A 2-tangle in LG(K) is such a subgraph. The particular case of cliques and anticliques subgraphs, with respect to common neighborhood inside the subgraph, is the case of digonal torsads $t^2$, which are particular 2-tangles. We can then generalize the 2-tangle to any type of graph:

**Definition:** *A 2-tangle of a graph G is an induced subgraph of G, where the set of edges with one end vertex owned by the 2-tangle and the other end owned by the complementary set of vertices in G, is a bipartite complete subgraph of G.*

This definition is stable in N(G). It means that, if t is a 2-tangle of G and T is a 2-tangle of N(G), there is an isomorphism mapping 2-tangles from G and N(G) where T(N(G)) = N(t(G)). We could write T=N(t).

**III - Realizability and construction**

**III - 1 Realizability of Gauss code and Gauss diagram**

Any knot diagram G(K) of a knot on the sphere S² is represented by a crossing sequence which is a 2-occurrence word W(K), a chord diagram CD(K), or an enlacement graph LG(K). It would be more exact at this stage to denote these representations W(G(K)), CD(G(K)), LG(G(K)). As any attribute can be chosen as an arbitrary decoration, the study of realizability can be restricted to a shadow of the knot, then to lacets or alternating knot diagrams.

Nevertheless many 2-occurrence words, chord diagrams and graphs are not the representations of knots. It is necessary to characterize these graphs. In the earliest time of the knot theory Gauss defined the chord diagram and found a necessary condition (i) equivalent to the statement that LG(K) is Eulerian Later on in 1936, Dehn found a sufficient algorithmic solution based on the existence of a touch Jordan curve which is the image of a transformation of the knot diagram by successive splits replacing all the crossings [De]. A long time after in 1976, Lovasz and Marx found a second necessary condition (ii) and finally during the same year, Rosenstiehl found the third condition (iii) which allowed the set of these three conditions to be sufficient. The last characterization is based on the tripartition of graphs into cycles, cocycles and bicycles. The algebraic concept covers and demonstrates the three conditions (i), (ii), (iii) as the orthogonality conditions for cycles and cocycles but even more by an extension made in 2001 by Rosenstiehl and Crapo which extended the law to general conditions for lacets on any type of 2-manifolds (orientable or not) when the complement of the closed curve is a 2-colorable family of discs. The last condition defines the lacet as a particular embedded curve.

It is remarkable to notice that the characterization of a 2-occurrence word and of a chord diagram is made by a characterization of their enlacement (intersection) graph. Denote by N(x) the set of neighbours of the vertex x in LG(K).

**Theorem 2-1 :** *A sequence W is a valid crossing sequence on the 2-sphere or on the Euclidean plane if and only if LG(W) is such that:*
*(i) For every vertex x, N(x) is even (Gauss).*
*(ii) For any pair (x,y) which is not an edge of LG(W), N(x) ∩ N(y) is even (Lovasz & Marx).*
*(iii) The edges (x,y) for which N(x) ∩ N(y) is even form a cocycle of LG(W) (Rosenstiehl).*

The lacet, Eulerian circuit of G(K), can be built as a diagonal curve (Eulerian line graph) of the checkerboard plane graph (or of its dual). The faces of G(K) are 2-colored and we remind that the checkerboard graph has one vertex on each face of one color and an edge is joining two incident faces of the same color. Another checkerboard graph is built under the same rule but with the other color and each graph is the dual of the other.

**Theorem 2-2 :** *If W is a valid crossing sequence with r connected components then there exist on the sphere $2^{(r-1)}$ distinct curves having W as crossing sequence. (Rosenstiehl)*

Connected components are corresponding to connected sets of chords of a CD and disjoint subgraphs of a LG. For each component there is the choice between one of the two checkerboard graphs. See figure 9. One chord diagram can then have several lacet representations when the number of components is not 1, i.e. when the knot for which the lacet is the shadow is not prime. Nevertheless the colors of the faces of G(K) are exchanged when we rotate the knot component 180° around an axis joining the 2 segments of strand attaching the component to the remaining part. This rotation transforms checkerboard component to its dual but, as this rotation is an ambient isotopy the resulting knot diagram is of the same knot type.

**Theorem 2-3 :** *A chord diagram has representations by lacets which are shadows of the same knot.*

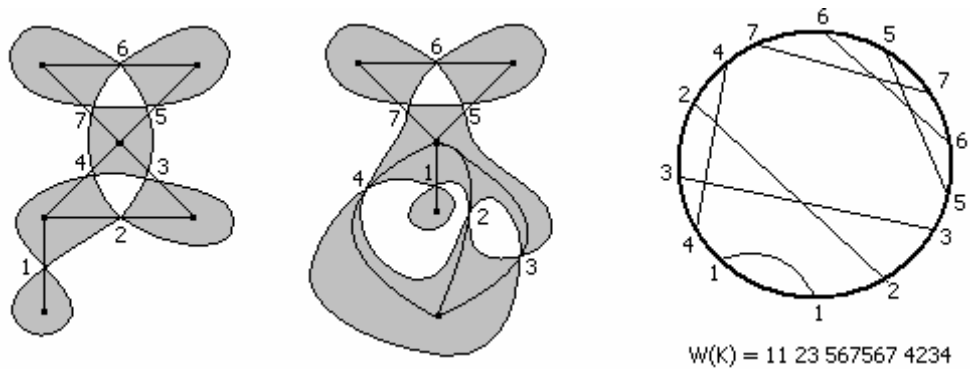

Figure 9 . Different lacets with their checkerboard graphs (one of the two dual versions of the whole and two of the four dual versions of the subgraph of two components) corresponding to the same chord diagram of a 3-component knot. 2 lacets are represented among 4 possible ones on the sphere.

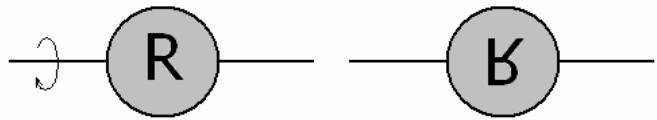

Figure 10 . Flip of a component : the colors of the faces of the component are exchanged and the checkerboard subgraph of the component is transformed to its dual.

On the knot diagram G(K), the flip is a flip of the corresponding subgraph, leading to an isomorphism of G(K), a different embedding of the abstract planar graph. Here is an example where a knot move by ambient isotopy is equivalent to an isomorphism of its diagram.

**III-2  Construction of a lacet or a knot diagram from theGauss code**

P. Rosenstiehl and R.C.Read exhibit a natural algorithm to draw a lacet by following the crossing sequence of the Gauss code with choice of turning loops to the left or to the right. They express its limitation as the number of possibilities grows exponentially like $2^{(n/3)}$ [ReRo].

Another algorithm is given by Dehn splitting method leading us to a contact curve. It is shown that any contact curve can be represented by a planar bipartite CD. The crossing curve is built from the plane representation of this CD. This algorithm is $O(n^2)$.

Another algorithm makes use of a pile of twin stacks and is $O(n)$. It was published in 1982 by P. Rosenstiehl and R.E. Tarjan [RoTa].

Above algorithms are not only valid for construction but also for testing realizability.

## IV - Flype and weighted chord diagrams

### IV-1 Definition

A flype is a 180° rotation of a 2-tangle as in figure 13.

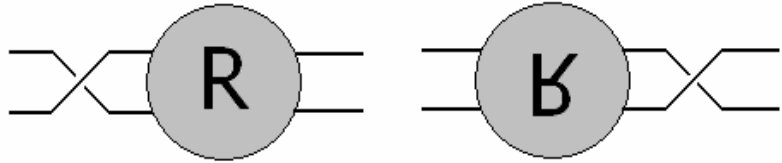

Figure 13 . Flype move

In the example of figure 13, we may consider that the disc supporting the tangle is rotating with it. The disc is cut from the surface of embedding and then glued back by surgery after rotation. Then, if all the over (resp. under) strands of the tangle are changed to under (resp. over) by the flype, the signs of the crossings who are representing the co-orientations with the disc remain unchanged. The axis of the rotation separates the 4 strands of the tangle into two sets of two. The "*active*" crossing is shifted unchanged from one side of the flype to the other and keep the same sign because the invariance of tangle-connectors by flype makes orientation of the strand invariant.

**Lemma 3-1 :** *Local writhe is invariant by flype*.

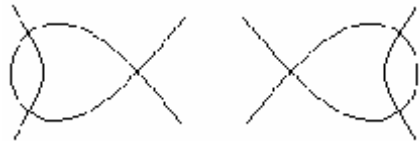

Figure 14 . The smallest non trivial flype is the "hornfish" where the tangle is a simple digon

Where a flype is possible we call it an *opportunity of a flype (OF)*. For each OF there is a 2-strand cycle of OFs built with active crossings and active 2-tangles. In particular, each 2-tangle has a complementary 2-tangle for which we will extract the active crossings like in figure 14.

On figure 15, four successive flypes are processed up to be back to the initial position. By a closer look one could notice that the third position is equivalent to the first one by S-isotopy, and also the fourth is S-isotopic to the second one. The number of steps for a cycle of flypes is even in any case and half of the positions are S-isotopic of the other half.

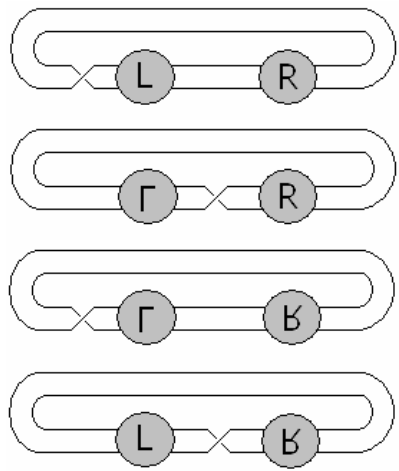

Figure 15 . A cycle of flypes

One of the oldest conjectures of knot theory was posed by P.G. Tait in 1898 and was proved by W. Menasco and M. Thistlethwaite in 1993 [MeTh]. Tait conjecture for alternating knots states that:

*Two alternating knots are isotopic if and only if any two corresponding diagrams on $S^2$ are related by a finite sequence of flypes.*

This process generates at most a finite number of projections. If a representation of a knot distinguishes all different types of knots and is invariant by any flype, then it is a complete invariant of alternating knots and that is the goal of this article. So let us consider how flypes act on chord diagrams and enlacement graphs.

**IV-2  Flype on the chord diagram**

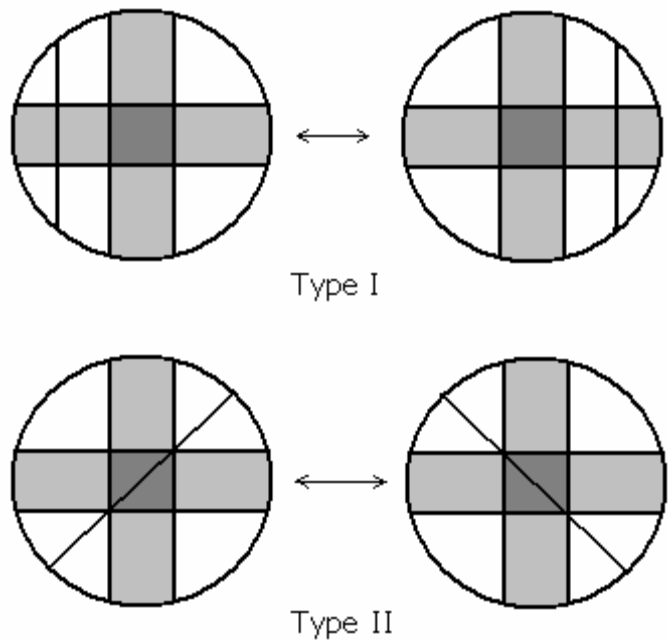

Figure 16 . The two types of flypes

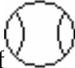

Figure 17 . Diagrams of the two types of flypes . Left and right
pairs of ends are joined .

Type I configuration is endowing exactly one tangle of 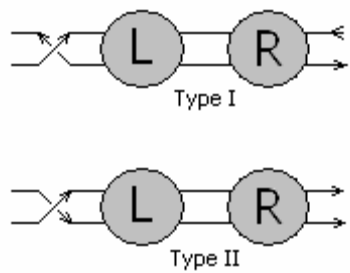 connector class. More than one would disconnect the lacet into several components. Less than one is impossible to obtain the right orientation of the strings. On a CD, this particular tangle is secant to all the others which are parallel. An active chord moving by flype through n tangles on G(K) has n different positions in both G(K) and CD(K), whatever is the type of the flype.

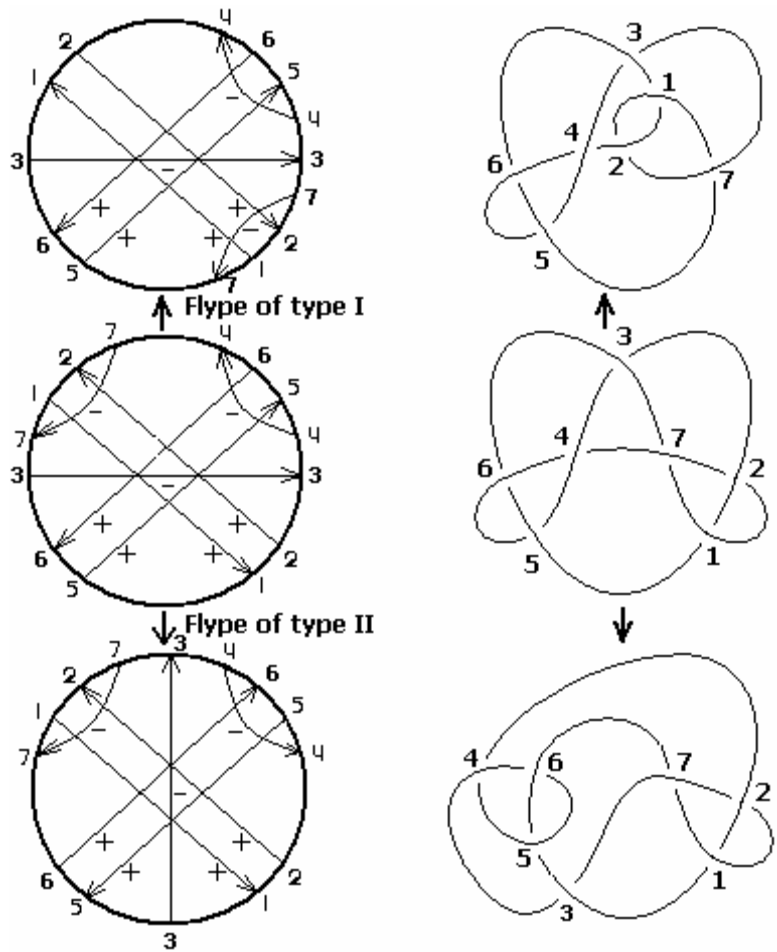

Figure 18 . Two types of flypes on $7_7$* knot

Two different types of knots must have different knot diagrams on S² and then different chord diagrams. The two conjugate mirror images of a knot diagram are distinguished by the crossing signs which are opposite on every crossing when orientation of the strand is kept. Then two chiral knots have the same (set of) CD but two different sign attributes. One alternating knot has different CDs, each one being the image of the other by a serie of flypes under the rule of Tait's conjecture. A signed CD distinguishes an alternating knot but is not invariant. Nevertheless we already proved that a CD is invariant by S-isotopy and component flipping. Flype is the only move of an alternating knot on the sphere. So, the only change of a signed CD is caused by the flype.

Now let's consider the geometrical state sum of all signed CDs of the same alternating knot. Each summand differ only by the position of active crossings (active chords). This state sum CD can then be represented by a chord-weighted CD where the active chord is duplicated in all possible positions by flype and where the weight of a chord is one if the chord is not active and the weight is 1/n if the chord is active and has n different positions by flype. For convenience of the notation we may use the sign of the crossing for the sign of the weight. The state sum signed CD is invariant by flype and distinguishes knots with chirality. Then it is a complete invariant representation of alternating knots.

**Theorem 3-2 :** *The chord-weighted chord diagram where*
*i ) sign of the weight of the chord is the local writhe of the corresponding crossing.*
*ii ) the absolute value of the weight is 1/n if and only if the chord has n different positions by flype.*
*iii ) every position of a chord by flype is represented by a chord*
*is a complete invariant graph representation of an alternating knot.*

We implied that a chord-weighted CD is a CD which is not rooted, not oriented and not labeled. For example see figure 18 where the chord-weighted chord diagram of $7_7$* is drawn with the three flype cycles, two of type I and one of type II. The weight is represented by a fraction or by a color.

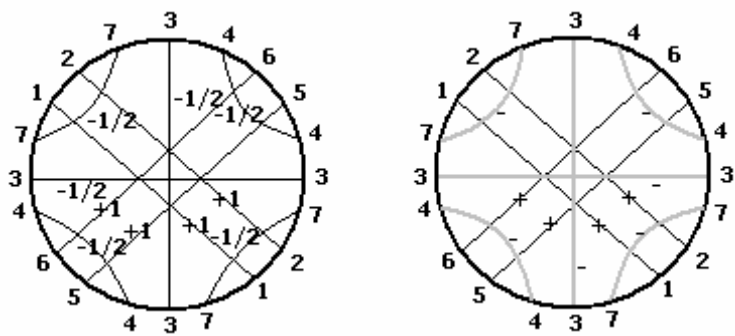

Figure 18 . Chord-weighted chord diagram of $7_7$* knot

Proof :

We have shown that the chord diagram is bijective with the shadow of a one-component knot on the 2-sphere. We have shown that the chord-weighted chord diagram (cwCD) is flype-invariant by construction of the weighted chords and by lemma 3-1.

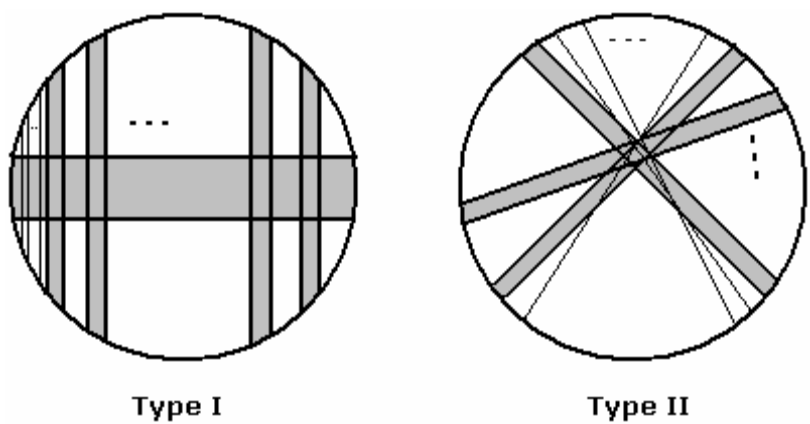

Figure 19 . p chords on q intervals between parallel (type I) or secant (type II) tangles

On figure 19 we generalize to p active chords on q intervals. For type I and for type II we have q tangles. On each interval of the cwCD we place one chord with +/- p/q weight, the sign being the sign of the crossing. Any non active chord has +/- 1 for weight. Nevertheless some active chords may also be represented by weighted chords of value +/-1 when p=q.

A condensed cwCD is also possible by attaching two signed attributes to the chords of the diagram, one indicating the sign and order of the $t^2$ torsad, the other indicating the sign of the crossing (the same along a torsad) and the averaged presence of the chord.

**IV-3  Chirality**

Local writhe (crossing sign) is inherited from the oriented knot diagram and assigned to cwCD as the sign of the weight of the chord. When mirroring a knot, the signs of all the crossings are inverted without to change anything else on cwCD. If the (unlabelled) cwCD is unchanged by the inversion of all the signs, then the knot is achiral, otherwise it is chiral. The test of chirality is made easy by the uniqueness of cwCD representation. Let's use three examples borrowed to C. Cerf Atlas [Ce] which are not easy to detect chirality by classical methods and just look at the figure 20. The inversion of the signs of the chords which corresponds to the mirroring is the identity for cwCD($8_{12}$) because of its symmetry. $8_{12}$ is then achiral. Mirroring on $10_{71}$ is changing its cwCD and then we can conclude easily that $10_{71}$ is chiral although the only method found until now to detect its chirality is complex. As simple is $8_{17}$ which appears to be amphichiral: mirror image is identity plus inversion of the orientation of the knot. $8_{17}$ is not invertible. For the same reason of symmetry and uniqueness of cwCD, we can also detect non-invertibility.

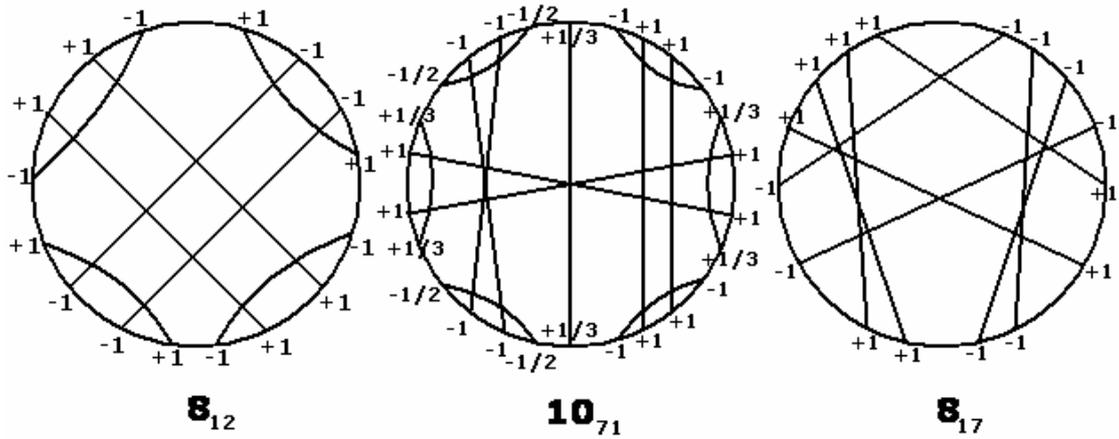

Figure 20 . cwCD($8_{12}$), cwCD($10_{71}$), cwCD($8_{17}$)

**IV-4 Symmetry**

For the same reason of uniqueness of cwCD representation, the symmetry of the knot is the one of cwCD and then easy to get.

**IV-5 Structure of flypes and 2-tangle decomposition of alternating knots**

Let B be a 3-ball and $t=t_1 \cup t_2$ a union of disjoint 2 arcs properly embedded in B. Then we call the pair (B, t) a *2-tangle*. Let S be a 2-sphere intersecting a knot K in 4 points. Then K is decomposed by S into two 2-tangles ($B_1$, $t_1$) and ($B_2$, t2), and the union ($B_1$, $t_1$) $\cup$ ($B_2$, t2) is called a *2-tangle decomposition* of K. If we consider now the 4-regular simple multigraph G(K) which is the knot diagram of K, a 2-tangle is an induced subgraph of G(K) of valence 4. When a closed curve is cutting the projection of the knot curve in 4 points and decomposing it into two 2-tangles, it is separating G(K) into two induced subgraphs of valence 4.

The simplest 2-tangles are :

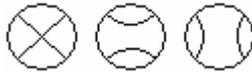

Any 2-tangle is seen from its exterior as one of these three trivial 2-tangles representing three classes of equivalence. We call them *connectors*. The only effect of a tangle on its complementary tangle in the knot is the orientation of the strand, and this orientation can be fixed by the connector of the tangle. It is easy to check that it is a relation of equivalence and each class of 2-tangle by this relation will be represented by its connector. A non-trivial 2-tangle must own at least one crossing in the first class on the list above and two crossings in the second and the third classes.

The structure of the composition of 2-tangles in a knot can be represented by a link if we replace each 2-tangle by a single crossing. An alternative and better structure is obtained by replacing each tangle by its connector. Then the nature of the knot is protected: number of components and orientation of each strands linking the tangles. The structure is a Contact-Crossing-Link. It is a recursive structure.

Flype opportunity (OF) is characterized by a sequence of a torsad $t^2_n$ and m 2-tangles and that is a $t^2_{n+m}$ torsad in the first tangle structure. Its closure is a (p, 2) torus link. The smallest OF is a sequence of one crossing and two 2-tangles where the (first) structure is the trefoil (closing of $t^2_3$).

The sum of all the tangles in the OF-cycle out of the active crossings is equivalent to:

flype of type I, for any number of active crossings,

flype of type II, for an odd number of active crossings,

flype of type II, for an even number of active crossings.

Other configurations are splitting the knot into 2 components. More generally, for n-tangles, the connector-class is a chord diagram with n chords. The connector will represent the action of a tangle on its complement to the knot.

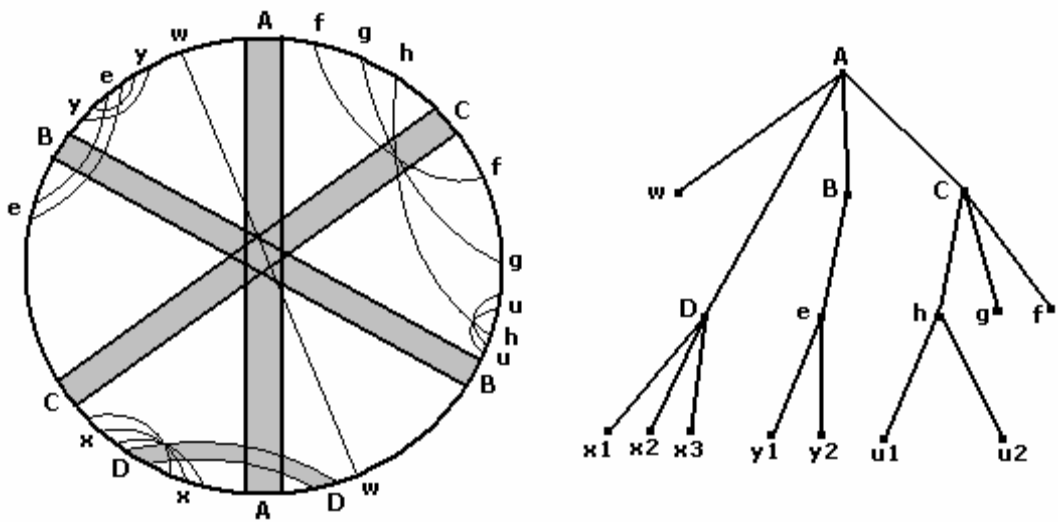

Figure 21 . Tangle tree

Tangles are drawn like thick chords on CD. We may find a tree structure to the OF-cycles where active crossings are on the leaves of the tree. No chord can be active on two different OF-cycles which are located on the leaves of the tree. See figure 21 where both type I and type II are represented at different levels (some flypes are trivial). The tree is not unique and is one of the spanning trees of the chord-intersection graph LG(K).

**V - Enlacement Graph and Mutations**

**V-1 Invariance and uncompleteness of LG**

LG is the intersection graph of chords on a CD. The representation of flype move has been described on the chord diagram of a knot on chapter IV-2 and it is clear on figure 16 that the two types of flype moves don't change the relations of intersections between chords. The number of chords is unchanged too and then vertices and edges of LG are unchanged by flypes of type I and II.

**Theorem 3-3 :** *The enlacement graph of a knot is invariant by flype move.*

Tait-Thistlesthwaite-Menasco theorem [Ta], [MeTh], shows that any knot diagram of the same alternating knot can be obtained from any other by a finite sequence of flype moves. Then:

**Corollary 3-4 :** *The enlacement graph of a knot is an invariant of an alternating knot.*

Nevertheless LG is not complete because different cwCDs and then different alternating knots may have the same enlacement graph. It is well known on Vassiliev-Kontsevich theory that the intersection graph conjecture is false (see [CDL] and [BMe] for more details) and it has been pointed out that mutant knots are not distinguished by this graph. Some types of LG graphs although were listed and proved to be in accordance with this conjecture. Here, mutation invariance of the chord-intersection graph is proven by the combinatorial geometry of knot diagrams and not by the 1-term and 4-term relations. Nevertheless the characterization of LG-equivalence is not well understood. We will bring an answer to the case of alternating knot isotopes. Some applications to Vassiliev-Kontsevich theory will be in a next print from the same author.

The enlacement graph doesn't detect mirroring without to add e local writhe attributes to vertices which are images of the chords. We'll not consider chirality in the sequel until we build a final graph based on some enhancement of LG.

The three symmetries acting on a tangle to obtain a mutant don't change the intersections of this tangle with other chords. Let's start with the classical example of Kinoshita-Terasaka and Conway on figure 22 below. We modified this pair of knots by changing some over-under attributes to make them alternating and then stay in our area of study. The last two extreme chords in the chain of torsads are flypable and are represented by dotted lines. The symmetry of the rotating tangle allows getting two mutants only instead of the maximum of four. That is the reason why we built a new set of mutants to put on evidence these four different mutants. It is reported on figure 23 after.

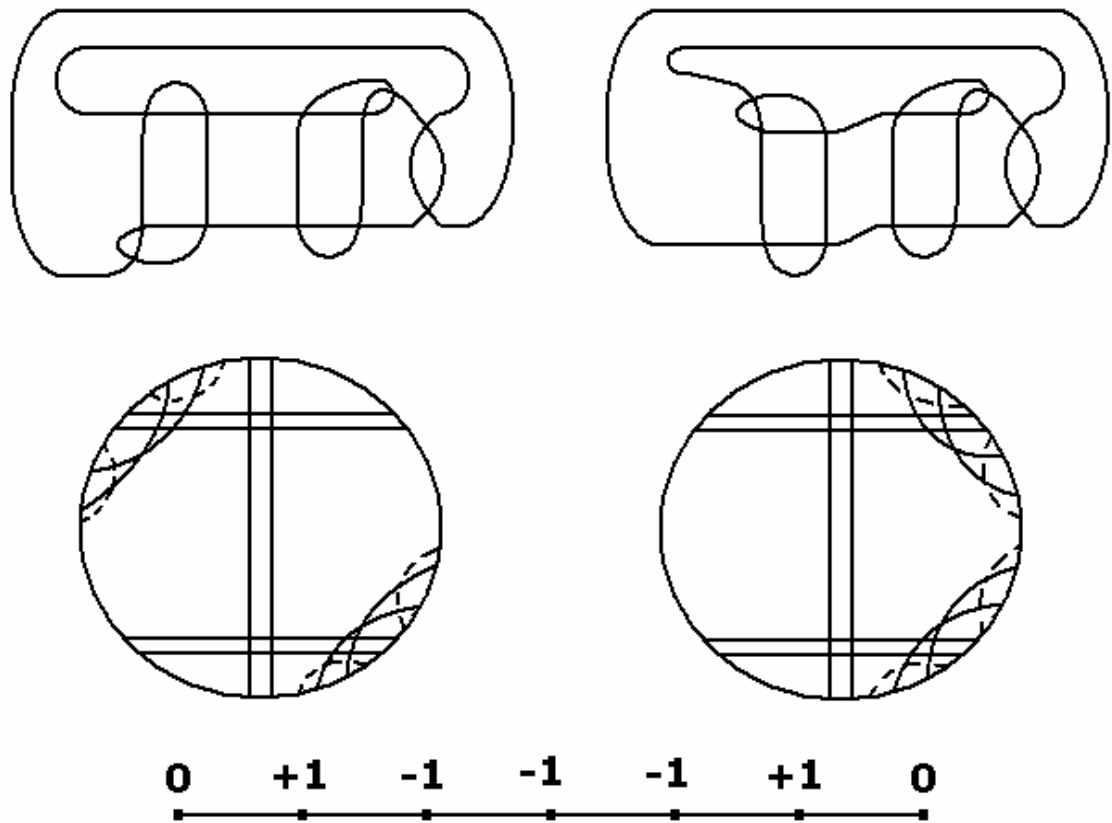

Figure 22 . Diagrams of one alternating variation of Kinoshita -Terasaka and Conway knots, chord-weighted cwCD chord diagrams, and common condensed enlacement graph.

Classical mutation is made by applying one of 3 symmetries to a 2-tangle of the knot, 4 if we encounter identity:
- Ÿ   Id : Identity
- Ÿ   H : Symmetry on an horizontal axis and over-under exchange on each crossing of the tangle
- Ÿ   V : Symmetry on a vertical axis and over-under exchange
- Ÿ   $\pi$ : Rotation counterclockwise of 180°

These four symmetries are the only ones which maintain the knot as a one component alternating knot. For example, a rotation of 90° split the knot into a link of 2 components and an axis symmetry without over-under switching produces a non-alternating knot. To be effective and not trivial, the 2-tangle and (with) its 2-tangle complement must not have the same symmetry as applied by the mutation. Such tangle internal symmetry includes the added symmetry of all flypes and we must consider the cwCD. For avoiding this difficulty and for getting a knot with 4 mutants, we chose our tangles by 2 asymmetric cuts on two different non algebraic knots, $8_{16}$ and $8_{17}$, which are rooted by the Borromean link. Then, to get a knot and not a link, we placed the 2 tangles with different connector-equivalence, here:

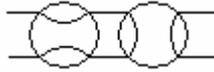

In fact, we know and it is easy to check that a closed sequence of two one-component 2-tangles is a knot if and only if their respective connectors are not the same.

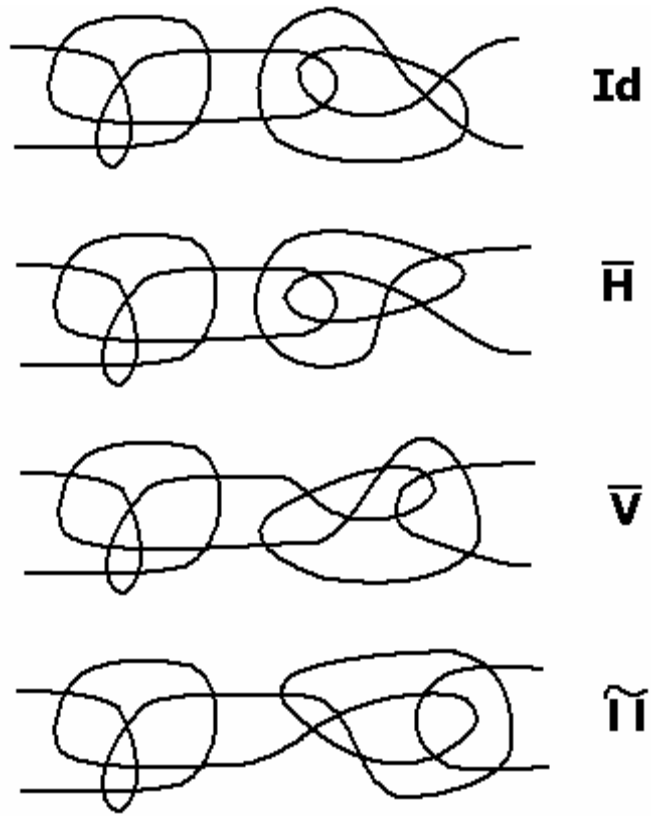

Figure 23 . Four mutants obtained from four symmetries on the same tangle.
Closing ends are not drawn.

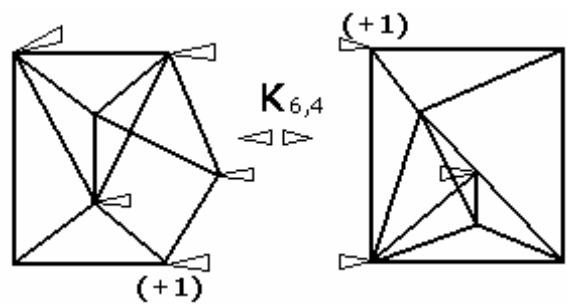

Figure 24 . Condensed enlacement graph of the 4 knots of figure 23
The bipartite complete bundle $K_{6,4}$ joining the 2 tangles is represented by small triangles.

## V-2 Enhancing the enlacement graph for mutant detection

### V-2-1 Local writhe enhancement

We may notice that the 3 symmetries of mutation leave invariant the orientation of the strands, the black-white partition of the faces of G(K). Moreover the co-orientation given by the ε attribute (local writhe on a crossing) is unchanged by axis-symmetry and rotation but ε sign is inverted if we switch over-under attribute.

$Id : \varepsilon \rightarrow \varepsilon$
$H : \varepsilon \rightarrow -\varepsilon$
$V : \varepsilon \rightarrow -\varepsilon$
$\pi : \varepsilon \rightarrow \varepsilon$

**Theorem 3-5 :** ε *local writhe distinguish H- or V-mutants from Id- or π-mutants.*

If we assign ε attribute as a weight for the vertices of LG(K) or for its condensed form (all crossings of a torsad have the same ε), we obtain a new graph $LG^\varepsilon$ which distinguish these two types of mutants.

**Corollary 3-6 :** *The graph $LG^\varepsilon(K)$ obtained from the enlacement graph of a knot K by assigning local writhe to each vertex of LG(K) distinguish H- or V-mutants from Id- or π-mutants.*

We notice that it is the case for Kinoshita-Terasaka-Conway knot which is obtained by H-symmetry on a (Id-) 2-tangle.
We must now distinguish (Id, H) from (V, π) if we want to be complete with classical mutation.

### V-2-2 Orientation enhancement

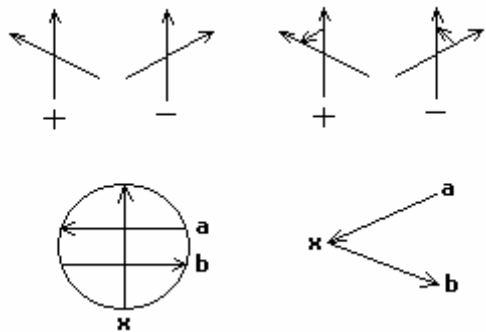

Figure 25 . Sign of a crossing path and orientation of chords on CD and of edges on LG

Let's give a sign to a crossing path like on figure 25. It is positive (resp. negative) sign when the second path is from right to left (resp. from left to right). This sign attribute has been used by several authors [CaEl] [DoTh] and give a more direct access to the embedding of the lacet. The two paths on one crossing have opposite signs and we draw an arrow from the positive strand to the negative strand. We do the same on the CD and the chord becomes to be an arrow. If an arrow chord intersect an other one from right to left (resp. from left to right), then the edge on LG will be directed from the first corresponding vertex to the second (resp. from the second to the first).

We consider now the direction of the edges of LG incident to a vertex of the tangle to rotate and a vertex of the complementary tangle. We note α this orientation and -α the opposite orientation.

$Id_{CD} : \alpha \rightarrow \alpha$
$H_{CD} : \alpha \rightarrow -\alpha$
$V_{CD} : \alpha \rightarrow \alpha$
$\pi_{CD} : \alpha \rightarrow -\alpha$

α-orientation distinguishes ($Id_{CD}$, $V_{CD}$) from ($H_{CD}$, $\pi_{CD}$)
A symmetry of mutation on CD(K) is not the same symmetry on G(K) but there is for each tangle a one to one correspondence between the two sets of symmetries ($Id_{CD}$, $H_{CD}$, $V_{CD}$, $\pi_{CD}$) and (Id, H, V, π), the last set operating on G(K).

### V-2-3 Complete invariant graph for alternating mutant knots

By merging the benefits of both enhanced LG graphs as described above, we can obtain a graph which distinguishes mutants.

ε local writhe distinguish H- or V-mutants from Id- or π-mutants.

(Id, π) = ($Id_{CD}$, $\pi_{CD}$) or ($\pi_{CD}$, $Id_{CD}$). Then {Id, π} = {$Id_{CD}$, $\pi_{CD}$}.
(H,V) = ($H_{CD}$, $V_{CD}$) or ($V_{CD}$, $H_{CD}$). Then {H,V}= {$H_{CD}$, $V_{CD}$}.
Then ε distinguish {$Id_{CD}$, $\pi_{CD}$} from {$H_{CD}$, $V_{CD}$}.
α-orientation distinguishes ($Id_{CD}$, $V_{CD}$) from ($H_{CD}$, $\pi_{CD}$).
Then (ε, α) distinguish $Id_{CD}$, $H_{CD}$, $V_{CD}$, $\pi_{CD}$ and due to bijection between the two symmetries sets, (ε, α) distinguish Id, H,V, π.

**Theorem 3-7 :** *The graph $LG^{\varepsilon,\alpha}(K)$ obtained from the enlacement graph of a knot K by assigning local writhe to each vertex of LG(K) and orienting its edges by the sign of crossing path, distinguishes all mutant alternating knots .*

### V-3  LG-equivalence characterization

Classical mutation as defined above is not the only knot change keeping the enlacement graph invariant. We may notice that a change on a set of chords in CD can be decomposed into a sequence of moves of connected sets of chords in correspondence with an induced subgraph of LG(CD). So we will restrict in the sequel to the relative position of connected sets of chords in a CD.

To keep condensed LG invariant, common neighborhoods must be kept unchanged whatever are the new intersections between chords. It is conjectured that only two types of intersections are made possible and it is partitioning the subgraph into two parts, the first being adjacent to the complementary part of the graph by a bundle of edges forming a complete bipartite subgraph, and the other being adjacent by no edges (a bundle forming an anti- complete bipartite subgraph). Such a subgraph is a 2-tangle. This definition of a 2-tangle is independent to its representation in a chord diagram. It is even well defined for any type of graph. We can then generalize tangles to abstract graphs.

Now consider that a 2-tangle is a subgraph of valence 4 in G(K). Any 2-tangle decomposition of G(K) is then supported by a link diagram G(TK) where every crossing point is assigned to a 2-tangle of the decomposition, including simple crossing points also as a special 2-tangle. Some crossing points of G(TK) have common neighborhoods and we can bundle them as in the above study into positive and negative torsads. If we permute crossing points inside such bundles, CD(TK) and LG(TK) are trivially unchanged. On the same manner, if we permute 2-tangles inside the same bundles, we don't change LG(K) but we change CD(K) and then the knot image of the CD. One particular move is a flype on TK which we call a *tangle flype*. In effect the active crossing point of a flype is replaced here by an active tangle. As any composition of permutations is a composition of transpositions, it is also a composition of tangle flypes.  If tangle flype looks like a flype on CD(TK), it is not the same on a diagram of TK where, due to the three different connectors, some crossing points are replaced by contact points. Then the diagram of TK becomes a contact-cross curve JK (J like Jordan) representing CD(JK) which is not realizable by a pure crossing curve. A new generalization of flype is then applied to contact points.

Let's consider applying tangle flypes step by step. The active tangle and its complementary have different connectors because it is a necessary and sufficient condition to form a knot and not a multi-component link. This property on the two connectors is maintained and so is the property to be a knot. Then tangle flype is well defined and the composition of tangle flypes (permutation of flypes on a bundle) is well defined also.

We call *LG-equivalence* of two representations of a knot or a lacet (K, G(K), CD(K),...) the relation for these representations to have the same enlacement graph LG.

**Definition :** *An extended mutation is an automorphism from the set of the knots to itself obtained by composition of :*
*1) classical mutations : symmetries on tangles by Id, H, V, or π*
*2) permutations of tangles in the same bundle (tangle flypes)*

A *bundle* is a torsad on KT diagram, which is a link representing the structure of a decomposition of the knot into tangles. The bundle is represented by secant or parallel chords in CD(JK).

**Theorem 3-8 :** *Two alternating knots have the same enlacement graph if they are isomorphic by extended mutation.*

**Corollary 3-9 :** *Extended mutants are LG-equivalent alternating knots.*

**Conjecture 3-10 :** *Two alternating knots have the same enlacement graph if and only if they are isomorphic by extended mutation.*

**Conjecture Corollary 3-11 :** *LG-equivalent alternating knots are extended mutants*

Let's illustrate by some examples and some representations.

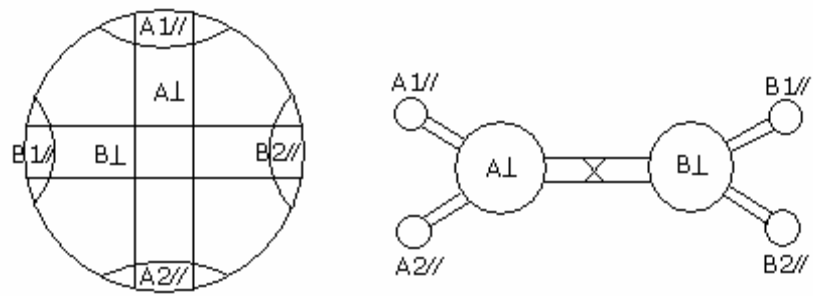

Figure 26 . The cross of tangles A and B on CD and on LG. A tangle may have up to 2 "bells". The part of the tangle with ends of each chord on 2 different pieces of strand is noted orthogonal. The parts of the tangle with ends of each chord on the same piece of strand are noted parallel (the bells). The complete bipartite subgraph bundling the 2 orthogonal parts of the tangles is noted with a cross.

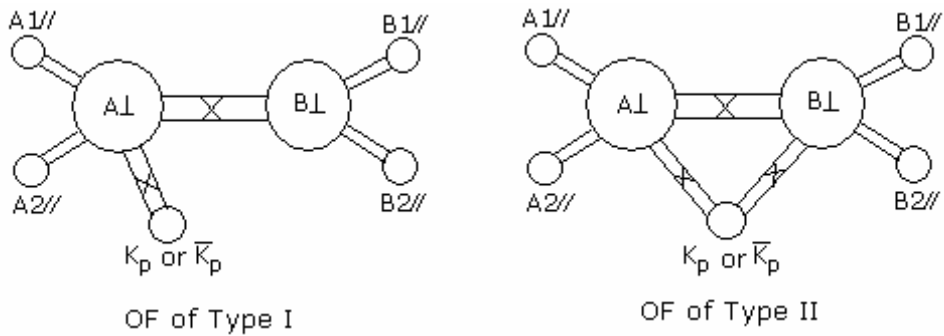

Figure 27 . The two types of flype opportunities on LG with active cliques or anticliques of order p.

The classical mutation may be considered to act on a tangle but it is equivalent to consider that it acts on the complementary tangle. We may also consider that it acts on the complete bipartite subgraph of LG bundling two tangles each one being complementary of the other. Then the four states obtained by applying Id, H, V, and π can be used like attributes of bipartite complete bundles on the enlacement graph.
The permutation of tangles on a chord diagram is changing the order of the ends of the tangles on the circle of the CD. We may apply this circular ordering to the corresponding LG. The problem is to give priority of ordering to the left or to the right of the circle. This ordering is not defined without to give orientation to the tangles or to the chords. We must then firstly give orientation to the chords like it is done on V-2-2. Then the abstract LG graph is becoming an embedded oriented graph when adding circular ordering to its edges.

By enriching a condensed LG with $Z_4$-weights on complete bipartite subgraphs and with orientation and circular ordering of the edges, we obtain an invariant complete graph of an alternating knot up to chirality. We must add e attribute (local writhe) on the vertices for detecting chirality.

Alternatively we can add circular ordering to the $Z_2$-weighted directed graph issued from LG (see paragraph V-2-3) and we get again an invariant complete graph of an alternating knot.

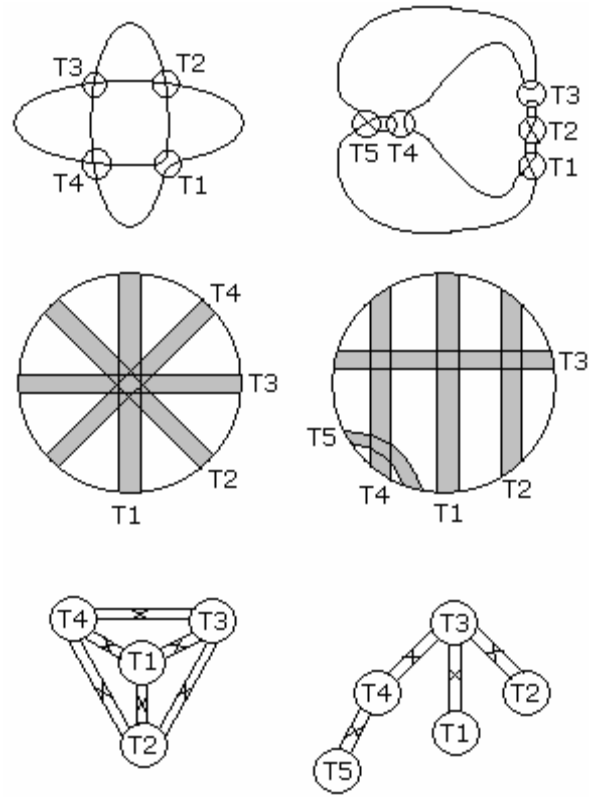

Figure 28 . Links of tangles TK and JK curves with CD(JK) and LG(JK). Tangles are represented on the diagram by their connector class. For example, a transposition between T1 and T2 with same connector position on the cycle of tangles, is an extended mutation (tangle flype) which doesn't change the intersections between chords and then leave LG invariant.

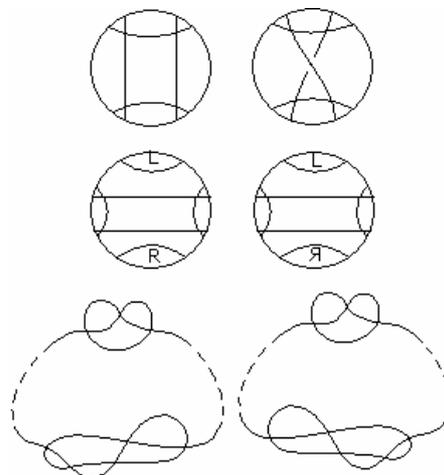

Figure 29 . Twisted mutation is possible for composite knots.

We pushed out of classical mutations the twist of the tangle on a CD because such an operation changes the intersection graph (edge complementation on the orthogonal part of the tangle). Nevertheless we can operate like this without to change LG if the tangle has an empty orthogonal part, i.e. if the tangle is made of two disjoint bells which means that these bells represent component knots. The twist mutation is inverting the Gauss word of one component like it is shown on figure 29.

## VI - Comments

Alternating knots are on a frontier between graph and knot theories and both advances on these two fields can contribute. Finding a convenient mathematical object which is a complete invariant of links and maybe of more general concepts would be the first step before to dive into any calculus, but the straight path cannot exist before to understand the deep nature of links through intense and diversified research. When limiting the problem to the simplest knots it was then proved that this way is not hopeless. Nevertheless, non-alternating knots have large amplitude of change of their shadow and of their related graphs. Another way must be found. The number of components of a link is of global nature where many actual theories are based on skein relations who are local and mix knots and links. Here chord diagrams can help to define the right algebra.

The author is although preparing a sequel on non-alternating knots and Vassiliev-Kontsevich theory.

It is worth to notice that some statistic laws can be represented by a state sum on a graph, like an integral of all possible states, but, unlike integration, it may be possible, as pointed out in our example, to retrieve all possible states from such an unique representation. We may conjecture that such state sum graphs could be defined sometimes in statistical mechanics, quantum physics, or other areas of science. A set of undefined states of parameters specifying a phenomenon, can be represented by one object only and processed like that. Then the uncertainty of the states is not an obstruction to process the phenomenon.

Christian Soulie is independent researcher. Any remark on this preprint is welcome at Chsoulie@aol.com or by mail at:
Christian Soulié
10, Allée du Parc
91 700 Fleury-Mérogis
France